\documentclass[12pt]{article}
\usepackage{amsmath,amsthm,amssymb
}

\allowdisplaybreaks

\newtheorem{theorem}{Theorem}[section]
\newtheorem{corollary}{Corollary}[section]

\theoremstyle{remark}

\theoremstyle{remark}

\theoremstyle{remark}
\newtheorem{remark}{Remark}[section]

\newcommand{\di}{\partial}
\newcommand{\C}{{\Bbb C}}
\newcommand{\N}{{\Bbb N}}
\begin{document}

\begin{center}{\Large \bf
   Lowering and raising operators for the free Meixner class of orthogonal polynomials}\end{center}

{\large Eugene Lytvynov}\\ Department of Mathematics,
Swansea University, Singleton Park, Swansea SA2 8PP, U.K.\\
e-mail: \texttt{e.lytvynov@swansea.ac.uk}\vspace{2mm}

{\large Irina Rodionova}\\
Department of Mathematics,
Swansea University, Singleton Park, Swansea SA2 8PP, U.K.\\
e-mail: \texttt{i.rodionova@swansea.ac.uk}\vspace{2mm}

{\small

\begin{center}
{\bf Abstract}
\end{center}

\noindent  We compare some properties of the lowering and raising operators for the classical and free classes of Meixner polynomials on the real line. } \vspace{2mm}

\vspace{1.5mm}

  \vspace{1.5mm}

\section{Classical and free Meixner classes}

In this note, we will compare some properties of the lowering and raising operators for the classical and free classes of Meixner polynomials on the real line.

Recall that the classical Meixner class  consists of all monic orthogonal polynomials $(P_n(x))_{n=0}^\infty$
whose exponential generating function has the form
\begin{equation}\label{treser}
\sum_{n=0}^\infty\frac{z^n}{n!}\,P_n(x)=\exp(x\Psi(z))f(z)=:G(x,z),
\end{equation}
where $z$ is from a neighborhood of zero in $\mathbb C$, $f$ and $\Psi$ are analytic functions in a neighborhood of zero such that $f(0)=1$, $\Psi(0)=\Psi'(0)=0$, and the measure of orthogonality of $(P_n(x))_{n=0}^\infty$, denoted by $\mu$, has  infinite support \cite{Meixner} (see also \cite{Chihara}).

Recall that, given a finite set $A$, a set partition of $A$ is a collection $\{A_1,\dots,A_m\}$ ($m\in\mathbb N$) of non-empty, mutually disjoint subsets of $A$ whose union is $A$. We denote by $\mathcal P_n$ the collection of all set partitions of the set $\{1,2,\dots,n\}$.
For any random variables $X_1, \dots, X_n$ on a probability space $(\Omega, \mathcal{A}, P)$ which have all their moments finite, the cumulant of $X_1, \dots, X_n$, denoted by $C_n(X_1, \dots, X_n)$, is defined recurrently through the formula \begin{equation}\label{skejgh}\mathbb E(X_1 \dotsm X_n)=\sum_{\pi\in\mathcal P_n}\prod_{A\in \pi}C(A,X_1,\dots, X_n),\end{equation} where for any $A=\{i_1,\dots,i_k\}\subset\{1,\dots,n\}$
\begin{equation}\label{ftufyytfry}C(A,X_1,\dots, X_n):=C_k(X_{i_1},\dots, X_{i_k}).\end{equation}
Then the cumulant generating function of a random variable $X$ is defined by
\begin{equation}\label{fytfdydt}
C_X(z)=\sum_{n=1}^\infty \frac{z^n}{n!}\, C_n(X,\dots,X),\end{equation}
where $z$ is from a neighborhood of zero in $\mathbb C$ for which the series in \eqref{fytfdydt} converges absolutely.
For a probability measure $\rho$ on $(\mathbb R,\mathcal B(\mathbb R))$ whose Laplace transform is analytic in a neighborhood of zero, the cumulant transform  of $\rho$, denoted by $C_\rho(z)$, is defined as $C_X(z)$, where the random variable $X$ has distribution $\rho$.
In fact, we have:
$$ C_\rho(z)=\log\bigg(\int_{\mathbb R}e^{zx}\,\rho(dx)\bigg).$$

The above assumptions on the polynomials from the Meixner class automatically imply that the generating function  $G(x,z)$ can be represented as
\begin{equation}\label{ersa}
G(x,z)=\exp(x\Psi(z)-C_\mu(\Psi(z))).
\end{equation}

As shown in  \cite{Meixner}, a system of orthogonal polynomials $(P_n(x))_{n=0}^\infty$ belongs to the Meixner class if and only if there exist $l\in\mathbb R$, $\lambda\in\mathbb R$, $t>0$, and $\eta\ge0$ such that the $(P_n(x))_{n=0}^\infty$ satisfy the recurrence relation
\begin{equation}\label{saera}
xP_n(x)=P_{n+1}(x)+(\lambda n-l)P_n(x)+n(t+\eta(n-1))P_{n-1}(x),\quad n\in\mathbb N_0:=\mathbb N\cup\{0\},\end{equation}
where $P_0(x)=1$  and $P_{-1}(x)=0$. It is easy to see that, if we set $l=0$ in \eqref{saera}, the corresponding measure $\mu$ will become centered, whereas $l\ne0$ corresponds to the shift of the centered measure by $l$. Therefore, we will restrict our attention to  centered measures, by setting $l=0$.

One has to distinguish the five following cases:

a) $\lambda=\eta=0$, $\mu$ is  Gaussian measure and $(P_n(x))_{n=0}^\infty$ are Hermite polynomials;

b) $\lambda\neq0$ and $\eta=0$, $\mu$ is  centered Poisson measure,  $(P_n(x))_{n=0}^\infty$ are Charlier polynomials;

c) $\eta>0$ and $\lambda^2=4\eta$, $\mu$ is  centered gamma measure,
$(P_n(x))_{n=0}^\infty$ are Laguerre polynomials;

d)  $\eta>0$ and $\lambda^2>4\eta$, $\mu$ is  centered Pascal measure, $(P_n(x))_{n=0}^\infty$ are Meixner polynomials of the first kind;

e)  $\eta>0$ and
$\lambda^2<4\eta$, $\mu$ is  Meixner measure,  $(P_n(x))_{n=0}^\infty$ are Meixner polynomials of the second kind, or Meixner--Polaczek polynomials in other terms.

In what follows, we will denote by $\mu_{\lambda,\eta}^{(t)}$, $G_{\lambda,\eta}^{(t)}$, $\Psi_{\lambda,\eta}^{(t)}$, and $C_{\lambda,\eta}^{(t)}$ the measure $\mu$ and the functions $G$, $\Psi$, and $C_\mu$, respectively, which correspond to the choice of the parameters $\lambda$, $\eta$, and $t$ as in \eqref{saera}. For $t=1$, we will usually skip the upper index ${}^{(1)}$.

For any possible choice of $\lambda$ and $\eta$, and for any $t>0$, we have $\Psi_{\lambda,\eta}^{(t)}(z)=\Psi_{\lambda,\eta}(z)$, $C_{\lambda,\eta}^{(t)}(z)=tC_{\lambda,\eta}(z)$, so that
\begin{equation}\label{dhltsrk}
G_{\lambda,\eta}^{(t)}(x,z)=\exp(x\Psi_{\lambda,\eta}(z)-tC_{\lambda,\eta}(\Psi_{\lambda,\eta}(z))).
\end{equation}

Furthermore, for any $\lambda$ and $\eta$ as above,
there exists a L\'evy process $X_{\lambda,\eta}=(X_{\lambda,\eta}(t))_{t\ge0}$ such that $X_{\lambda,\eta}(0)=0$ a.s., and for each $t>0$ the random variable $X_{\lambda,\eta}(t)$ has distribution $\mu_{\lambda,\eta}^{(t)}$  (see e.g.\ \cite{S}).

Let $\nu_{\lambda,\eta}$ be the probability measure on $\mathbb R$ which is the measure of orthogonality of the polynomials $(Q_n(x))_{n=0}^\infty$ satisfying
\begin{equation}
xQ_n(x)=Q_{n+1}(x)+\lambda(n+1)Q_n(x)+\eta n(n+1)Q_{n-1}(x),\quad n\in\mathbb N_0,
\label{dsar}\end{equation}
with $Q_0(x)=1$, $Q_{-1}(x)=0$.

Thus, if $\eta=0$, then $\nu_{\lambda,0}=\delta_\lambda$ (the Dirac measure with mass at $\lambda$), while for $\eta>0$ $\nu_{\lambda,\eta}$ is the measure of orthogonality of polynomials from the Meixner class satisfying \eqref{saera} with $\lambda$, $\eta$ as above,  $l=-\lambda$ and $t=2\eta$. Then, for $z$ from a neighborhood of zero in $\mathbb C$, we have:
\begin{align}
C_{\lambda,\eta}(z)&=\int_{\mathbb R}(e^{sz}-1-sz)s^{-2}\,\nu_{\lambda,\eta}(ds)\notag\\
&=\int_{\mathbb R}\sum_{n=2}^\infty \frac{s^{n-2}z^n}{n!}\,\nu_{\lambda,\eta}(ds),\label{ufytt}
\end{align}
i.e., $s^{-2}\,\nu_{\lambda,\eta}(ds)$ is the L\'evy measure of the L\'evy process $X_{\lambda,\eta}$.

In free probability, Meixner's class  of orthogonal polynomials was introduced and studied by Anshelevich  \cite{a1} (see also \cite{a2}) and Saitoh, Yoshida \cite{SY}. (In fact, such polynomials had already occurred in many places in the literature even before \cite{a1,SY}, see \cite[p.~62]{BB} and \cite[p.~864]{a2} for  bibliographical references.)  A deep study of multivariate free Meixner polynomials of non-commutative variables has been carried out by Anshelevich, see \cite{a3,a2,a5} and the references therein. Bo\.zejko and Bryc \cite{BB}  gave a characterization of free Meixner polynomials in terms of a regression problem.
We also refer to \cite{a1,a3,BW} for $q$-interpolation of the classical and free Meixner classes.

Below, in the free case, for many objects we will use  the same notations as those used for their counterparts in the classical case. However, it should always be clear from the context which objects are being meant.

The free Meixner class  consists of all monic orthogonal polynomials $(P_n(x))_{n=0}^\infty$
whose (usual) generating function has the form
\begin{equation}\label{dhlk}
\sum_{n=0}^\infty z^n\,P_n(x)=(1-x\Psi(z))^{-1}f(z)=:G(x,z),
\end{equation}
where $z$, $\Psi$ and $f$ satisfy the same assumptions as in the classical case, and the measure of orthogonality $\mu$ has  infinite support.
Recall the following notation from $q$-analysis: for each $q\in[-1,1]$, we define $[0]_q:=0$ and $[n]_q:=1+q+q^2+\dots+q^{n-1}$, for $n\in\mathbb N$, and $[n]_q!:=[1]_q[2]_q\dotsm[n]_q$. In particular, the free probability corresponds to $q=0$, in which case $[0]_0=0$ and $[n]_0=1$ for $n\in\mathbb N$. Thus the free analog of $\exp(x)=\sum_{n=0}^\infty \frac{x^n}{n!}$
is the resolvent function $(1-x)^{-1}=\sum_{n=0}^\infty x^n$, so that \eqref{dhlk} is indeed the free analog of \eqref{treser}.

Consider a non-commutative probability space which is a von Neumann algebra $\mathcal A$ with a normal, faithful, tracial state $\tau:\mathcal A\to \mathbb C$, i.e., $\tau(\cdot)$ is linear, continuous in
weak$*$ topology, $\tau(ab)=\tau(ba)$, $\tau(\operatorname{id})=1$, $\tau(aa^*)\ge0$, and $\tau(aa^*)=0$ implies $a=0$. A (non-commutative) random variable $X$ is a self-adjoint element of $\mathcal A$.

We  denote by $\mathcal {NC}(n)$ the collection of all non-crossing partitions of $\{1,\dots,n\}$, i.e., all set partitions $\pi=\{A_1,\dots,A_k\}$, $k\ge1$, of $\{1,\dots,n\}$ such that there do not exist $A_i,A_j\in\pi$, $A_i\ne A_j$, for which the following inequalities hold: $x_1<y_1<x_2<y_2$
for some $x_1,x_2\in A_i$ and $y_1,y_2\in A_j$. For any non-commutative random variables $X_1,\dots,X_n$, the free cumulant $C_n(X_1,\dots,X_n)$ is defined through formulas \eqref{skejgh}, \eqref{ftufyytfry} in which $\mathbb E$ is replaced by $\tau$, and $\mathcal P_n$ is replaced by
$\mathcal {NC}(n)$. Then, analogously to \eqref{fytfdydt}, the free cumulant generating function of a random variable $X$ is defined by $$C_X(z)=\sum_{n=1}^\infty z^n  C_n(X,\dots,X).$$
If $\rho$ is a probability measure on $(\mathbb R,\mathcal B(\mathbb R))$ with compact support, we define the free cumulant transform  $C_\rho(z)$ of $\rho$ as $C_X(z)$, where a random variable $X$ has distribution $\rho$, i.e., $\tau(X^n)=\int_{\mathbb R}x^n\,\rho(dx)$.

By \cite{a1}, the generating function of the fee Meixner polynomials can be represented as
\begin{equation}\label{gdxsghshtr}
G(x,z)=(1-x\Psi(z)+C_\mu(\Psi(z)))^{-1}
\end{equation}
(compare with \eqref{ersa}). Furthermore, a system of orthogonal polynomials $(P_n(x))_{n=0}^\infty$ belongs to the free Meixner class if and only if there exist $l\in\mathbb R$, $\lambda\in\mathbb R$, $t>0$, and $\eta\ge0$ such that the $(P_n(x))_{n=0}^\infty$ satisfy the recurrence relation
\begin{equation}\label{ewaew}
xP_n(x)=P_{n+1}(x)+(\lambda [n]_0-l)P_n(x)+[n]_0(t+\eta[n-1]_0)P_{n-1}(x),\end{equation}
and we again set $l=0$ in order to set the measure $\mu$ centered. Analogously to the classical case, we thus again have five classes of free Meixner polynomials.  Next, analogously to \eqref{dhltsrk}, we have
\begin{equation}\label{tsre}
G_{\lambda,\eta}^{(t)}(x,z)=(1-x\Psi_{\lambda,\eta}(z)+tC_{\lambda,\eta}(\Psi_{\lambda,\eta}(z)))^{-1}
\end{equation}
(we have used obvious notations).

Recall that non-commutative random variables $X_1,\dots,X_n$ are called freely independent if, for any $k\ge1$ and any non-constant choice of $Y_1,\dots,Y_k\in\{X_1,\dots,X_n\}$, $C_k(Y_1,\dots,Y_k)=0$ \cite{Speicher} (see also \cite{BB}). A non-commutative stochastic process $X=(X(t))_{t\ge0}$ is called a free
L\'evy process if $X_0=0$ and  the increments of $X$ are freely independent and stationary, see \cite{biane} for details.

By \cite{a1,SY}, for any allowed parameters $\lambda$ and $\eta$,
there exists a free  L\'evy process $X_{\lambda,\eta}=(X_{\lambda,\eta}(t))_{t\ge0}$ such that,
 for each $t>0$, the random variable $X_{\lambda,\eta}(t)$ has distribution $\mu_{\lambda,\eta}^{(t)}$.
Let $\nu_{\lambda,\eta}$ be the probability measure on $\mathbb R$ which is the measure of orthogonality of the polynomials $(Q_n(x))_{n=0}^\infty$ satisfying
\begin{align}
xQ_n(x)&=Q_{n+1}(x)+\lambda[n+1]_0Q_n(x)+\eta [n]_0[n+1]_0Q_{n-1}(x),\notag\\
&=Q_{n+1}(x)+\lambda Q_n(x)+\eta Q_{n-1}(x),
\label{dghcdfr}\end{align}
with $Q_0(x)=1$, $Q_{-1}(x)=0$. Thus, if $\eta=0$, then $\nu_{\lambda,0}=\delta_\lambda$, while for $\eta>0$ 
\begin{equation}\label{gutgu}
\nu_{\lambda,\eta}=\mu_{\lambda,0}^{(\eta)}.\end{equation}
 (Thus, for $\eta>0$, $\nu_{\lambda,\eta}$ is the free Gaussian distribution if $\lambda=0$ and the free Poisson distribution if $\lambda\ne0$.)
Then, for $z$ from a neighborhood of zero in $\mathbb C$, we have:
\begin{align}
C_{\lambda,\eta}(z)&=\int_{\mathbb R}((1-sz)^{-1}-1-sz)s^{-2}\,\nu_{\lambda,\eta}(ds)\notag\\
&=\int_{\mathbb R}\sum_{n=2}^\infty z^ns^{n-2}\,\nu_{\lambda,\eta}(ds),\label{hyfytds}
\end{align}
i.e., $s^{-2}\,\nu_{\lambda,\eta}(ds)$ is the free L\'evy measure of the free L\'evy process $X_{\lambda,\eta}$ \cite{bnt3}.

\begin{remark} The polynomials of the classical Mexiner class with $\eta>0$ naturally appear  in the study of  a realization of the renormalized  square of white noise, see \cite{AFS}. In fact, as it (indirectly) follows from \cite{sn} (see also \cite{BL}),  the polynomials of the free Meixner class with $\eta>0$ are related to the renormalized square of free white noise. 
\end{remark}

\section{Raising and lowering operators}

We fix any $\lambda\in\mathbb R$ and $\eta\ge0$. In the case of classical Meixner polynomials, we define the raising and lowering operators, $\di_{\lambda,\eta}^\dag$ and $\di_{\lambda,\eta}$, as linear operators given through
$$  (\di_{\lambda,\eta}^\dag P_n)(x)=P_{n+1}(x),\quad (\di_{\lambda,\eta} P_n)(x)=nP_{n-1}(x),$$
where $(P_n(x))_{n=0}^\infty$ satisfy \eqref{saera} (with $t=1$ and $l=0$). We will denote by $x$ 
the operator of multiplication by variable  $x$ in $L^2(\mathbb R,\mu_{\lambda,\eta})$.
By \eqref{saera}, we then have 
\begin{equation}\label{huftyfi}
x=\di_{\lambda,\eta}^\dag (1+\lambda\di_{\lambda,\eta}+\eta\di_{\lambda,\eta}^2)+\di_{\lambda,\eta}.
\end{equation}
Since the operators $\di_{\lambda,\eta}^\dag$, $\di_{\lambda,\eta}$ are unbounded in $L^2(\mathbb R,\mu_{\lambda,\eta})$, we will later identify their domain.

In the free case, we define the raising and lowering operators, $\di_{\lambda,\eta}^\dag$ and $\di_{\lambda,\eta}$ through
$$  (\di_{\lambda,\eta}^\dag P_n)(x)=P_{n+1}(x),\quad (\di_{\lambda,\eta} P_n)(x)=[n]_0P_{n-1}(x),$$
where $(P_n(x))_{n=0}^\infty$ satisfy \eqref{ewaew} (again with $t=1$ and $l=0$). The operator of multiplication by $x$ then has the same representation \eqref{huftyfi} in $L^2(\mathbb R,\mu_{\lambda,\eta})$. Note that the operators $x$, $\di_{\lambda,\eta}^\dag$, and $\di_{\lambda,\eta}$ are now bounded.

\subsection{Lowering operator}

We start with the classical case.
So, denote by $\mathcal P$ the set of all polynomials on $\mathbb  R$. For each $f(x)=\sum_{n=0}^N f_n P_n(x)\in\mathcal P$, $n\in\mathbb N_0$, and each $q\in\mathbb N$, we set
$$\|f\|_q^2=\sum_{n=0}^\infty |f_n|^2 (n!)^2 2^{nq}.$$
Let $H_q$ denote the Hilbert space obtained as the completion of $\mathcal P$ in the $\|\cdot\|_q$ norm. We then define the nuclear space
$\Phi:=\projlim_{q\to\infty} H_q$.
Since $(P_n(x))_{n=0}^\infty$ is a Scheffer system of polynomials,  there exists $q\in\mathbb N$ such that $H_q$ (and so $\Phi$) is continuously embedded into $L^2(\mathbb R,\mu_{\lambda,\eta})$ (see \cite{KSS}).

Denote by ${\cal E}_{\mathrm {min}}^1(\C)$ the set of all
entire functions on $\C$  of first order of growth and of
minimal type, i.e., a function $f$ entire on $\C$
belongs to ${\cal E}_{\mathrm {min}}^1(\C)$ if and only if for each $\varepsilon >0$ there exists $C>0$ such that $|f(z)|\le C\exp(\varepsilon|z| )$.
 Denote by ${\cal E}_{\mathrm
{min}}^1(\mathbb R)$ the set of restrictions to $\mathbb R$ of
functions from ${\cal E}_{\mathrm{min}}^1(\C)$. Following
\cite{KSS, KSWY}, we then  introduce  a norm on ${\cal E}_{\mathrm
{min}}^1(\C)$, and hence on  ${\cal E}_{\mathrm {min}}^1(\mathbb R)$, as follows. Each $f\in
{\cal E}_{\mathrm {min}}^1(\C)$ can be uniquely represented
in the form $ f(z)=\sum_{n=0}^\infty  f_n z^
n$,  and we set, for any $q\in\N$, $$ N_{
q} (f):=\sum_{n=0}^\infty |f_n|^2(n!)^2\,
2^{nq}.$$ By \cite[Theorems~2.5, 3.8 and subsec.~6.2]{KSS},
$\Phi={\cal E}_{\mathrm {min}}^1(\mathbb R)$
and the
two systems of norms on $\Phi$:
$(\|\cdot\|_{q},\  q\in\mathbb N)$ and $(N_{q}(\cdot),\ q\in\mathbb N)$ are equivalent, and
hence determine the same topology on $\Phi$.

Using the $\|\cdot\|_q$ norms, we easily conclude that the operators $\di_{\lambda,\eta}^\dag$ and $\di_{\lambda,\eta}$ act continuously on ${\cal E}_{\mathrm {min}}^1(\mathbb R)$.
We have the following theorem which describes the action of $\di_{\lambda,\eta}$ in the classical case
(compare with \cite[Theorem~2.2]{R} and
\cite[Theorem~4.1]{Ly}):

\begin{theorem}\label{fytfytf}
In the classical case, we have, for any $f\in{\cal E}_{\mathrm {min}}^1(\mathbb R)${\rm :}
\begin{equation}\label{gftyk} (\di_{\lambda,\eta}f)(x)=\int_{\mathbb R}\frac{f(x+s)-f(x)}s\,\nu_{\lambda,\eta}(ds),\quad x\in\mathbb R.\end{equation}
\end{theorem}

\begin{remark}\label{uyfyrfuyk}{\rm In the Gaussian case, i.e., when $\lambda=\eta=0$, we have $\nu_{0,0}=\delta_0$ and formula
\eqref{gftyk}  is understood in the limiting sense: $(\di_{0,0}f)(x)=(Df)(x):=f'(x)$.
}\end{remark}

\begin{remark}\label{hugfk}{\rm In the Poisson case, i.e., when $\eta=0$ and $\lambda\ne0$, we have
$\nu_{\lambda,0}=\delta_\lambda$, so that
$$ (\di_{\lambda,0}f)(x)=\frac{f(x+\lambda)-f(x)}{\lambda}\, .$$
Therefore, if $\eta\ne0$,
$$ (\di_{\lambda,\eta}f)(x)=\int_{\mathbb R}(\di_{\lambda,0}f)(x)\nu_{\lambda,\eta}(d\lambda),$$
and so the operator  $\di_{\lambda,\eta}$ is the $\nu_{\lambda,\eta}$-mixture of  the operators $\di_{\lambda,0}$. }\end{remark}

\noindent{\it Proof of Theorem} \ref{fytfytf}.
For each $q\in\mathbb N$, denote by $G_q$ the Hilbert space obtained as the completion of
${\cal E}_{\mathrm {min}}^1(\mathbb R)$ in the $N_q(\cdot)$ norm.
 As easily seen, the operator $\di_{\lambda,\eta}$ acts continuously in each $H_q$.  Hence, there exist $q_1\ge q_2\ge1$ such that $G_{q_2}$ is continuously embedded into $L^2(\mathbb R,\mu_{\lambda,\eta})$ and $\di_{\lambda,\eta}: G_{q_1}\to G_{q_2}$ is a continuous operator.  Choose $\varepsilon>0$ such that, for each $z\in\mathbb R$, $|z|<\varepsilon$, $e^{x z}$ and $e^{x\Psi_{\lambda,\eta}(z)}$ belong to $G_{q_1}$ as functions of $x$. Hence, $G(\cdot,z)\in G_{q_1}$ and
$$ (\di_{\lambda,\eta}G_{\lambda,\eta}(\cdot,z))(x)=zG_{\lambda,\eta}(x,z),$$
which implies
\begin{equation}\label{ygufr} \di_{\lambda,\eta}e^{xz}=\Psi_{\lambda,\eta}^{-1}(z)e^{xz}.\end{equation}
On the other hand, by \eqref{ufytt},
\begin{align}\int_{\mathbb R}\frac{e^{(x+s)z}-e^{xz}}{s}\,\nu_{\lambda,\eta}(ds)&=e^{xz}\int_{\mathbb R}\frac{e^{sz}-1}{s}\,\nu_{\lambda,\eta}(ds)\notag\\
&=e^{xz}C'_{\lambda,\eta}(z).\label{igty}\end{align}
By \cite[Proposition 1]{a5} (see also \cite{Meixner}), we have
\begin{equation}\label{jhfy}\Psi_{\lambda,\eta}^{-1}(z)=C'_{\lambda,\eta}(z).\end{equation}
Hence, \eqref{ygufr}--\eqref{jhfy} imply that \eqref{gftyk} holds when $f(x)=e^{xz}$, $|z|<\varepsilon$.

Since $\int_{\mathbb R}e^{\varepsilon|s|}\nu_{\lambda,\eta}(ds)<\infty$, there exits $C\ge1$ such that, for all $n\in\mathbb N$,
$$
\int_{\mathbb R}|s|^n\nu_{\lambda,\eta}(ds)\le C^n n!\,.
$$
Hence, for each $f(x)=\sum_{n=0}^\infty f_n x^n\in {\cal E}_{\mathrm {min}}^1(\mathbb R)$ and $x\in\mathbb R$,
\begin{align}\int_{\mathbb R}\bigg|\frac{f(x+s)-f(x)}{s}\bigg|\nu_{\lambda,\eta}(ds)
&\le \sum_{n=1}^\infty |f_n|\int_{\mathbb R}
\bigg|\frac{(x+s)^n-x^n}{s}\bigg|\nu_{\lambda,\eta}(ds)\notag\\
&\le \sum_{n=1}^\infty |f_n| \sum_{k=0}^{n-1}\binom{n}{k}|x|^k\int_{\mathbb R}|s|^{n-k-1}\nu_{\lambda,\eta}(ds)\notag\\
&\le \sum_{n=1}^\infty |f_n| n!\, n C^n (1\vee|x|)^n\notag\\
&\le N_{q_3}(f)\bigg(\sum_{n=1}^\infty 4^{-n}\bigg)^{1/2}<\infty,\label{ugfytft}
\end{align}
where $q_3=q_3(x)\in\mathbb N$ is chosen so that $q_3\ge q_2$ and $(4C(1\vee|x|))^2\le2^{q_3}$. 

Now, let $\{f^{(k)}\}_{k=1}^\infty\subset G_{q_3}$ be such that each $f^{(k)}$ is a linear combination of functions $e^{xz}$ with $|z|<\varepsilon$ and  $f^{(k)}\to f$ as $k\to\infty$ in $G_{q_3}$. Then, analogously to \eqref{ugfytft}, we conclude that
$$ \int_{\mathbb R}\frac{f^{(k)}(x+s)-f^{(k)}(x)}{s}\,\nu_{\lambda,\eta}(ds)\to 
\int_{\mathbb R}\frac{f(x+s)-f(x)}{s}\,\nu_{\lambda,\eta}(ds), $$
and since $\di_{\lambda,\eta}f^{(k)}\to\di_{\lambda,\eta}f$ in $G_{q_1}$, $(\di_{\lambda,\eta}f^{(k)})(x)\to(\di_{\lambda,\eta}f)(x)$. From here the theorem follows. \quad $\square$

We will now derive a free counterpart of Theorem~\ref{fytfytf}. 
By the theory of Jacobi matrices (see e.g.\ \cite{Ber}), in the free  case, the measure $\mu_{\lambda,\eta}^{(t)}$ is concentrated on the interval $[-(\lambda\vee\sqrt{\eta+t}),(\lambda\vee\sqrt{\eta+t})]$. 
Below $C^1([a,b])$ denotes the set of all continuously differentiable functions on interval $[a,b]$.

\begin{theorem}\label{u8iutg8tg} In the free case, we have for any $f\in C^1([-r_{\lambda,\eta},r_{\lambda,\eta}])$,
\begin{align} (\di_{\lambda,\eta}f)(x)&=\int_{[-r_{\lambda,\eta},r_{\lambda,\eta}]}\frac{f(x)-f(s)}{x-s}\,\nu_{\lambda,\eta+1}(ds)\label{rdtdts}\\
&=\int_{[-r_{\lambda,\eta},r_{\lambda,\eta}]}\frac{f(x)-f(s)}{x-s}\,\mu_{\lambda,0}^{(\eta+1)}(ds)
,\quad x\in\mathbb R.\label{tdser}\end{align}
Here, $r_{\lambda,\eta}:=(\lambda\vee\sqrt{\eta+1})$. 
\end{theorem}

\begin{remark}\label{ghfyrd}{\rm 
Unlike in the classical case, the integral representation of the operator $\di_{\lambda,\eta}$ in the free case uses the measure $\nu_{\lambda,\eta+1}$ with the `shifted' parameter $\eta+1$. So, in particular, in the free Gaussian case ($\lambda=\eta=0$),  $\di_{0,0}$ is {\it not\/} the operator of free differentiation 
$f(x)\mapsto (D_{\text{free}}f)(x)=\frac{f(x)-f(0)}{x}$. In fact, in the free Gaussian and free Poisson cases, i.e., when $\eta=0$, the integration on the right hand side of formula 
\eqref{tdser} is with respect to the measure of orthogonality $\mu_{\lambda,0}$. Thus, in particular, 
$$ P_{n-1}(x)=\int_{[-r_{\lambda,0},r_{\lambda,0}]}\frac{P_n(x)-P_n(s)}{x-s}\,\mu_{\lambda,0}(ds),$$
where $(P_n(x))_{n=0}^\infty$ are orthogonal with respect to $\mu_{\lambda,0}$. 
}\end{remark}

\begin{remark}\label{gyyfr}{\rm 
In the classical and free cases, denote 
by $m_{\lambda,\eta}(n)$ the $n$-th moment of $\nu_{\lambda,\eta}$:
$$m_{\lambda,\eta}(n):=\int_{\mathbb R}s^n\nu_{\lambda,\eta}(ds).$$ Then, by Theorem~\ref{fytfytf},
$$ \di_{\lambda,\eta}x^n
=\sum_{k=0}^{n-1}\binom{n}{k}m_{\lambda,\eta}(n-1-k)x^k,$$
while in the fee case, by Theorem~\ref{u8iutg8tg},
$$ \di_{\lambda,\eta}x^n
=\sum_{k=0}^{n-1}m_{\lambda,\eta+1}(n-1-k)x^k,$$
Recall the the free ($q=0$) analog of the binomial coefficient $\binom nk$
is 1.
}\end{remark}

{\it Proof of Theorem\/} \ref{u8iutg8tg}.
First, we note that, by \eqref{gutgu}, 
formulas \eqref{rdtdts} and \eqref{tdser}
 are equivalent. 
By \cite{a1}, for sufficiently small $z$
\begin{align}
\Psi_{\lambda,\eta}(z)&=\frac{z}{1+\lambda z+\eta z^2},\label{guyfu}\\
C_{\lambda,\eta}(\Psi_{\lambda,\eta}(z))&=
\frac{z^2}{1+\lambda z+\eta z^2}\, .\label{reazrwea}
\end{align}
Hence, by \eqref{gdxsghshtr}, we easily see that 
\begin{align*}
G_{\lambda,\eta}(x,z)&=\bigg(1-\frac{xz-z^2}{1+\lambda z+\eta z^2}\bigg)^{-1}\\
&=\frac{1+\lambda z+\eta z^2}{1+\lambda z+(\eta+1) z^2}\,(1-x\Psi_{\lambda,\eta+1}(z))^{-1}.
\end{align*}
Since
$$ (\di_{\lambda,\eta}G_{\lambda,\eta}(\cdot,z))(x)=zG_{\lambda,\eta}(x,z),$$
we, therefore, have
\begin{equation}\label{sres} \di_{\lambda,\eta}(1-xz)^{-1}=\Psi_{\lambda,\eta+1}^{-1}(z)(1-xz)^{-1}.
\end{equation} 
Next, by \eqref{hyfytds} 
\begin{align}
&\int_{[-r_{\lambda,\eta},r_{\lambda,\eta}]}\frac{(1-xz)^{-1}-(1-sz)^{-1}}{x-s}\,\nu_{\lambda,\eta+1}(ds)\notag\\
&\qquad =
(1-xz)^{-1}\int_{[-r_{\lambda,\eta},r_{\lambda,\eta}]}z(1-sz)^{-1}\nu_{\lambda,\eta+1}(ds)\notag\\
&\qquad=(1-xz)^{-1}z^{-1}\int_{[-r_{\lambda,\eta},r_{\lambda,\eta}]}\sum_{n=2}^\infty
z^ns^{n-2}\nu_{\lambda,\eta+1}(ds)\notag\\
&\qquad=(1-xz)^{-1}z^{-1} C_{\lambda,\eta+1}(z).\label{tyde}
\end{align}
By \cite[Proposition~1]{a5} (see also \cite{a1}),
\begin{equation}\label{xdd}\Psi^{-1}_{\lambda,\eta}(z)=
z^{-1}C_{\lambda,\eta}(z)
\end{equation}
(compare with \eqref{jhfy} and note that $z^{-1}C_{\lambda,\eta}(z)$ is the free derivative of $C_{\lambda,\eta}$). By \eqref{sres}--\eqref{xdd}, equality \eqref{rdtdts} holds for $f(x)=(1-xz)^{-1}$ for all sufficiently small $z$. From here, the general case follows by an easy approximation argument.\quad $\square$ 

\begin{remark}{\rm Note that, in view of \eqref{ygufr} and \eqref{sres}, in the classical case
$$\di_{\lambda,\eta}=\Psi_{\lambda,\eta}^{-1}(D),$$
 while in the free case
$$\di_{\lambda,\eta}=\Psi_{\lambda,\eta+1}^{-1}(D_{\text{free}}).$$
}\end{remark}

\subsection{Raising operator}

We again start with the classical case (cf.\ \cite{Ly,R}). Following \cite{Meixner}, we define  $\alpha,\beta\in\C$ through the
equation 
\begin{equation}\label{fytf}1+\lambda x+\eta
x^2=(1-\alpha x)(1-\beta x),\end{equation} or equivalently
$$\alpha+\beta=-\lambda, \quad\alpha\beta=\eta,$$
where in the case $\eta=0$ and $\lambda\neq0$ one sets
$\alpha=-\lambda$ and $\beta=0$. Evidently, the condition $\lambda\in\mathbb R$, $\eta\geq0$ is satisfied if
and only if either $\alpha,\beta\in\mathbb R$ and $\alpha$ and $\beta$ are
of the same sign, or $\operatorname{Im}(\alpha)\neq0$ and $\alpha$
and $\beta$ are complex conjugate.

\begin{theorem}\label{kjigyug}
In the classical  case, there exists $\varepsilon>0$ such that, for all $z\in\mathbb R$, $|z|<\varepsilon$, we have
\begin{align}
\di_{\lambda,\eta}^\dag e^{xz}&=
\left(\frac{x}{1+\lambda \Psi^{-1}_{\lambda,\eta}(z)+\eta(\Psi^{-1}_{\lambda,\eta}(z))^2}
-\frac{\Psi^{-1}_{\lambda,\eta}(z)}{1+\lambda \Psi^{-1}_{\lambda,\eta}(z)+\eta(\Psi^{-1}_{\lambda,\eta}(z))^2}
\right)e^{xz}\label{det6sw}\\
&=\bigg[x
e^{-z(\alpha-\beta)}\bigg(1+\frac{\alpha}{\alpha-\beta}\big(e^{z(\alpha-\beta)}-1\big)
\bigg)^2 \notag \\
&\quad-e^{-z(\alpha-\beta)}
\bigg(1+\frac{\alpha}{\alpha-\beta}\big(e^{z(\alpha-\beta)}-1\big)\bigg)
\frac{1}{\alpha-\beta}\big(e^{z(\alpha-\beta)}-1\big)\bigg]e^{xz},
\label{bvn7ed}\end{align} for $\alpha=\beta$ the formula being
understood in the limiting sense.
\end{theorem}

\begin{remark}
{\rm Theorem \ref{kjigyug} shows that the operator $\di_{\lambda,\eta}^\dag$ is represented through the operator of multiplication by $x$ and through analytic functions (in a neighborhood of zero) of the operator of differentiation $D$. 
}\end{remark}

{\it Proof of Theorem }\ref{kjigyug}.
By \eqref{huftyfi} and \eqref{ygufr}, we have
$$ xe^{xz}=(1+\lambda \Psi^{-1}_{\lambda,\eta}(z)+\eta(\Psi^{-1}_{\lambda,\eta}(z))^2)\di_{\lambda,\eta}^\dag e^{xz}+\Psi^{-1}_{\lambda,\eta}(z) e^{xz},$$
from where \eqref{det6sw} follows. To derive \eqref{bvn7ed} from \eqref{det6sw}, use \eqref{fytf} and the following formula (see \cite{Meixner}): 
$$ \Psi_{\lambda,\eta}^{-1}(z)=\bigg(\frac{1}{\alpha-\beta}\big(e^{z(\alpha-\beta)}-1
\big)\bigg)\bigg(1+
\frac{\alpha}{\alpha-\beta}\big(e^{z(\alpha-\beta)}-1\big)\bigg)^{-1}.\quad\square$$

Using formula \eqref{bvn7ed} and analogously to the proof of Theorem~\ref{fytfytf}, one easily derives  explicit formulas for the action of $\di_{\lambda,\eta}^\dag$ (cf.\ \cite{Ly,R}). Before formulating 
this result, we introduce the following natations: For each $s\in\mathbb R$, $s\ne0$ and for $f:\mathbb R\to\mathbb R$, we define
\begin{align*}
(\nabla_sf)(x):=&\frac{f(x+s)-f(x)}{s}\, ,\\
(U_sf)(x):=&f(x+s).
\end{align*}
Clearly, $D$, $\nabla_s$, and $U_s$ act continuously on ${\cal E}_{\mathrm {min}}^1(\mathbb R)$.

\begin{corollary}\label{ftsw5t}
We have the following representation of the 
 operator $\di_{\lambda,\eta}$ on ${\cal E}_{\mathrm {min}}^1(\mathbb R)$: 
for  $\lambda=\eta=0$
$$\di_{0,0}^\dag= x- D,$$
for  $\lambda\ne0$ and $\eta=0$
$$\di_{\lambda,0}^\dag= x(1-\lambda\nabla_{\lambda})-\nabla_\lambda,$$
for  $\eta>0$ and $\lambda^2=4\eta$
$$ \di_{\lambda,\eta}^\dag=x(D-1)^2-D(D-1),$$
and for $\eta>0$ and $\lambda^2\ne4\eta$
$$ \di_{\lambda,\eta}^\dag=x(1+\alpha\nabla_{\alpha-\beta})^2U_{\beta-\alpha}-(1+\alpha\nabla_{\alpha-\beta})\nabla_{\alpha-\beta}U_{\beta-\alpha}.$$
\end{corollary}

We proceed to consider the free case. 

\begin{theorem}\label{ygfyufr}
In the free case, there exists $\varepsilon>0$ such that, for all $z\in\mathbb R$, $|z|<\varepsilon$, we have
\begin{align}&
\di^\dag_{\lambda,\eta}(1-xz)^{-1}\notag\\
&=
\left(\frac{x}{1+\lambda \Psi^{-1}_{\lambda,\eta+1}(z)+\eta(\Psi^{-1}_{\lambda,\eta+1}(z))^2}
-\frac{\Psi^{-1}_{\lambda,\eta+1}(z)}{1+\lambda \Psi^{-1}_{\lambda,\eta+1}(z)+\eta(\Psi^{-1}_{\lambda,\eta+1}(z))^2}
\right)(1-xz)^{-1}\label{tydfyd}\\
&=\left(
x\, \frac{4z^2(\eta+1)^2}
{\big(2\eta+1+\lambda z+\sqrt{(1-\lambda z)^2-4z^2(\eta+1)}\big)
\big(1-\lambda z-\sqrt{(1-\lambda z)^2-4z^2(\eta+1)}\big)}\right.\notag\\
&\left.\qquad\quad+\frac{2z(\eta+1)}{2\eta+1+\lambda z+\sqrt{(1-\lambda z)^2-4z^2(\eta+1)}}\right)(1-xz)^{-1}. 
\label{yfutfr}
\end{align}
\end{theorem}

\begin{remark}
{\rm By Theorem \ref{ygfyufr},  the operator $\di_{\lambda,\eta}^\dag$ is represented through the operator of multiplication by $x$ and through analytic functions (in a neighborhood of zero) of the operator of free differentiation $D_{\text{free}}$. 
}\end{remark}

{\it Proof of Theorem }\ref{ygfyufr}. 
The derivation of  \eqref{tydfyd} is analogous to the classical case. So, we only have to show that \eqref{yfutfr} holds. By \eqref{guyfu} and  \eqref{tydfyd},
\begin{equation}\label{utg8u}
\di_{\lambda,\eta}^\dag(1-xz)^{-1}=\left(x\,
\frac{\Psi_{\lambda,\eta}(\Psi^{-1}_{\lambda,\eta+1}(z))}{\Psi^{-1}_{\lambda,\eta+1}(z)}-\Psi_{\lambda,\eta}(\Psi^{-1}_{\lambda,\eta+1}(z))\right)(1-xz)^{-1}.
\end{equation}
 Next, by \eqref{guyfu}, 
\begin{align*} \frac1{\Psi_{\lambda,\eta+1}(z)}&=\frac{1+\lambda z+\eta z^2}z+z\\
&=\frac1{\Psi_{\lambda,\eta}(z)}+z.\end{align*}
Hence,
$$\frac1z=\frac{1}{\Psi_{\lambda,\eta+1}(\Psi^{-1}_{\lambda,\eta+1}(z))}=\frac1{\Psi_{\lambda,\eta}(\Psi_{\lambda,\eta+1}^{-1}(z))}+\Psi_{\lambda,\eta+1}^{-1}(z),$$
from where
\begin{equation}\label{rtdtrde}
\Psi_{\lambda,\eta}(\Psi_{\lambda,\eta+1}^{-1}(z))=
\frac1{\frac1z-\Psi_{\lambda,\eta+1}^{-1}(z)}.
\end{equation}
Since
\begin{equation}\label{uftfu}\Psi_{\lambda,\eta+1}^{-1}(z)=\frac{1-\lambda z-\sqrt{(1-\lambda z)^2-4z^2(\eta+1)}}{2z(\eta+1)},
\end{equation}
by \eqref{rtdtrde}
\begin{equation}\label{rdsrt}
\Psi_{\lambda,\eta}(\Psi_{\lambda,\eta+1}^{-1}(z))=
\frac{2z(\eta+1)}{2\eta+1+\lambda z+\sqrt{(1-\lambda z)^2-4z^2(\eta+1)}}\,. 
\end{equation}
Now, \eqref{yfutfr}  follows from \eqref{utg8u}, \eqref{uftfu}, and \eqref{rdsrt}.\quad $\square$

\begin{remark}
 In the free case, it is still an open problem whether one can derive any  explicit formulas for the action of $\di_{\lambda,\eta}^\dag$ on a general function $f$.
\end{remark}

\begin{remark} In \cite{BL}, a study of a free Meixner class of    orthogonal  polynomials of infinitely many non-commutative variables has been initiated. We expect that the results of this note related to the free case may be generalized to this infinite dimensional setting, compare with \cite{Ly}. We hope that this will be discussed in \cite{BL2}.  

\end{remark}

{\bf Acknowledgements}. We would like to thank Marek Bo\.zejko for numerous useful discussions.  EL was partially supported by  an LMS Scheme 3 grant and by the PTDC/MAT/67965/2006 grant, University of Madeira.


\begin{thebibliography}{99}

\bibitem{AFS} L. Accardi, U. Franz, and M.  Skeide, Renormalized squares of white noise and other non-Gaussian noises as L\'evy processes on real Lie algebras, {\it Commun.\ Math.\ Phys.}\ {\bf 228} (2002)  123--150. 


\bibitem{a1}  M. Anshelevich, Free martingale polynomials,  {\it J. Funct.\ Anal.}\ {\bf 201} (2003) 228--261. 

\bibitem{a3} M. Anshelevich, Appell polynomials and their relatives, {\it  Int.\ Math.\ Res.\ Not.}\ {\bf 2004},
no.~65, 3469--3531. 


\bibitem{a2} M. Anshelevich, Free Meixner states,  {\it Commun.\ Math.\ Phys.}\ {\bf 276} (2007)  863--899. 


\bibitem{a5} M. Anshelevich, Orthogonal polynomials with a resolvent-type generating function, {\it Trans.\ Amer.\ Math.\ Soc.}\ {\bf 360} (2008)  4125--4143. 

\bibitem{bnt3}  O. E. Barndorff-Nielsen, S.  Thorbj\o rnsen, The L\'evy-It\^o decomposition in free probability, {\it  Probab.\ Theory Related Fields}\/ {\bf 131} (2005) 197--228.

\bibitem{Ber}  Ju.\ M. Berezanskii,  {\it Expansions in Eigenfunctions of Selfadjoint Operators}   (American Mathematical Society, 1968).


\bibitem{biane}  P. Biane, Processes with free increments, {\it Math.\ Z.}\ {\bf 227} (1998)  143--174. 


\bibitem{BB} M. Bo\.zejko and W. Bryc, On a class of free L\'evy laws related to a regression problem,  {\it J. Funct.\ Anal.}\ {\bf 236} (2006) 59--77. 




\bibitem{BW} W. Bryc and J. Weso{\l}owski,  Conditional moments of $q$-Meixner processes, {\it Probab.\ Theory Related Fields} {\bf 131} (2005) 415--441. 

\bibitem{BL} M. Bo\.zejko and E. Lytvynov, Meixner class of  non-commutative generalized stochastic processes with freely independent values I. A characterization, Arxiv preprint, 2008. 

\bibitem{BL2} M. Bo\.zejko and E. Lytvynov, Meixner class of  non-commutative generalized stochastic processes with freely independent values II. The generating function, in preparation. 

\bibitem{Chihara} T. S. Chihara, {\it An Introduction to Orthogonal Polynomials}
  (Gordon and Breach Sci.\ Pbl.,  1978).
  
\bibitem{KSS} Y. Kondratiev, J. L. Silva, and L. Streit, Generelized Appell systems,
{\it Methods Funct.\ Anal.\ Topology\/} {\bf  3} (1997), no. 3,
28--61.

\bibitem{KSWY} Y. Kondratiev, L. Streit, W. Westerkamp, and J.
Yan,  Generalized functions in infinite dimensional analysis,
{\it Hiroshima Math. J.} {\bf 28} (1998) 213--260.

\bibitem{Ly} E. Lytvynov, Polynomials of Meixner's type in infinite
dimensions---Jacobi fields and orthogonality measures, {\it J.
Funct. Anal.}\ {\bf 200} (2003) 118--149.

\bibitem{Meixner} J. Meixner, Orthogonale Polynomsysteme mit einem besonderen
Gestalt der erzeugenden Funktion, {\it J. London Math.\ Soc.}\ {\bf
9} (1934) 6--13.

\bibitem{R} I. Rodionova, 
Analysis connected with generating functions of exponential type in one and infinite dimensions,
{\it Methods Funct.\ Anal.\ Topology} {\bf  11} (2005) 275--297.

\bibitem{SY}  N. Saitoh and H. Yoshida, The infinite divisibility and orthogonal polynomials with a constant recursion formula in free probability theory, {\it Probab.\ Math.\ Statist.}\ {\bf 21} (2001)  159--170.


\bibitem{S} W. Schoutens, {\it Stochastic Processes and Orthogonal
Polynomials}, Lecture Notes in Statist., Vol.~146
(Springer-Verlag,  2000).

\bibitem{sn} P. \'Sniady, Quadratic bosonic and free white noises,  {\it Commun. Math.\ Phys.}\ {\bf 211} 
(2000)   615--628.


\bibitem{Speicher}  R. Speicher, Free probability theory and non-crossing partitions,
 {\it S\'em. Lothar.\ Combin.}\  {\bf 39} (1997), Art. B39c, 38 pp.\   (electronic)



\end{thebibliography}
\end{document}